\documentclass[12pt,a4paper]{article}

\usepackage{indentfirst}

\usepackage{amsmath,amssymb,amsfonts,amsthm,graphics}
\usepackage{makeidx}
\begin{document}

\title{\bf \Large Upper bounds involving parameter $\sigma_2$\\
for the rainbow connection\footnote{Supported by NSFC.} }
\author{\small Jiuying Dong, Xueliang Li\\
\small Center for Combinatorics and LPMC-TJKLC \\
\small Nankai University, Tianjin 300071, China\\
\small Email: jiuyingdong@126.com; lxl@nankai.edu.cn}

\date{}

\maketitle{}

\begin{abstract}

For a graph $G$, we define $\sigma_2(G)=min \{d(u)+d(v)| u,v\in
V(G), uv\not\in E(G)\}$, or simply denoted by $\sigma_2$. A
edge-colored graph is rainbow edge-connected if any two vertices are
connected by a path whose edges have distinct colors, which was
introduced by Chartrand et al. The rainbow connection of a connected
graph $G$, denoted by $rc(G)$, is the smallest number of colors that
are needed in order to make $G$ rainbow edge-connected. We prove
that if $G$ is a connected graph of order $n$, then $rc(G)\leq
6\frac{n-2}{\sigma_2+2}+7$. Moreover, the bound is seen to be tight
up to additive factors by a construction mentioned by Caro et al. A
vertex-colored graph is rainbow vertex-connected if any two vertices
are connected by a path whose internal vertices have distinct
colors, which was recently introduced by Krivelevich and Yuster. The
rainbow vertex-connection of a connected graph $G$, denoted by
$rvc(G)$, is the smallest number of colors that are needed in order
to make $G$ rainbow vertex-connected. We prove that if $G$ is a
connected graph of order $n$, then $rvc(G)\leq
8\frac{n-2}{\sigma_2+2}+10 $ for $2\leq \sigma_2\leq 6, \sigma_2\geq
28 $, while for $ 7 \leq \sigma_2\leq 8, 16\leq \sigma_2\leq 27$, $
rvc(G)\leq \frac{10n-16}{\sigma_2+2}+10$, and for $9 \leq
\sigma_2\leq 15, rvc(G)\leq \frac{10n-16}{\sigma_2+2}+A(\sigma_2)$
where $ A(\sigma_2)= 63,41,27,20,16,13,11,$ respectively.\\[3mm]
\end{abstract}

\section{ Introduction}

All graphs in the paper are finite, undirected and simple. Let
$\sigma_2(G)=min \{d(u)+d(v)| u,v\in V(G), uv\not\in E(G)\}$, or
simply denoted by $\sigma_2$.  The distance between two vertices $u$
and $v $ in $ G $, denoted by $d(u, v)$, is the length of a shortest
path between them in $ G$. The eccentricity of a vertex $v$ is
$ecc(v) :=\max_{x\in V(G)}d(v, x)$. The diameter of $G$ is $diam(G)
:= \max_{x\in V(G)}ecc(x)$. For the notations and terminology not
defined here, we follow the book Bollob\'{a}s [2].

A path in an edge colored graph with no two edges sharing the same
color is called a rainbow path. An edge colored graph is said to be
rainbow connected if every pair of vertices is connected by at least
one rainbow path. Such a coloring is called a rainbow coloring of
the graph. The minimum number of colors required to rainbow color a
connected graph is called its rainbow connection number, denoted by
$rc(G)$. Note that disconnected graphs cannot be rainbow colored and
hence the rainbow connection number for them is left undefined. A
natural and interesting quantifiable way to strengthen the
connectivity requirement was introduced by Chartrand et al. in [6].
An easy observation is that if $G$ has $n$ vertices then $rc(G)\leq
n-1$. Also, clearly, $rc(G) \geq diam(G)$ where $diam(G)$ denotes
the diameter of $G$.

It was shown by Chakraborty et al. [4] that computing the rainbow
connection number of an arbitrary graph is NP-Hard. To rainbow color
a graph, it is enough to ensure that every edge of some spanning
tree in the graph gets a distinct color. There have been attempts to
find better upper bounds in terms of other graph parameters like
connectivity, minimum degree and radius etc. Caro et al. [3] have
proved that if $\delta \geq 3$ then $rc(G)=\alpha n$ where $\alpha <
1$ is a constant. They conjectured that $\alpha =\frac{3}{4}$
suffices and proved that $\alpha <\frac{5}{6}$. They also proved
$rc(G)\leq (ln\delta/\delta)n(1+o_\delta(1))$. Krivelevich and
Yuster [7] have obtained the best known bound of
$\frac{20n}{\delta}$ using a strengthened connected two-step
dominating set. Later, Chandran et al. [5] used a connected two-step
dominating set to show that for every connected graph on $n $
vertices with minimum degree $\delta$ the rainbow connection number
is upper bounded by $3n/(\delta + 1) + 3$. This solves an open
problem from Schiermeyer [8]. The result nearly settles the
investigation for an upper bound of rainbow connection number in
terms of minimum degree which was initiated by Caro et al. [3].

Since the parameter $\sigma_2(G)$ plays an extremely useful role in
the studying of graph connectivity, Hamiltonian property, etc, it is
interesting to use $\sigma_2$ to study the $rc(G)$ of a graph $G$.
We are encouraged and motivated by the above ideas, results and
proof methods. We give a upper bound of $rc(G)$ as a function of
$\sigma_2(G)$, which is stated as the following Theorem 1.

\vskip 0.2 cm \noindent{\bf Theorem 1.}  For a connected graph $G$
of order $n$, $rc(G)\leq 6\frac{n-2}{\sigma_2+2}+7$.
 \vskip 0.2 cm

The following examples show that our bound  $rc(G)\leq
6\frac{n-2}{\sigma_2+2}+7$ are seen to be tight up to additive
factors.

\noindent{\bf Example 1:} Add edges to $K_{2, \sigma_2/2-1}$ such
that the part $\overline{K}_{\sigma_2/2-1}$ of $K_{2, \sigma_2/2-1}$
is $K_{\sigma_2/2-1}$, we denote the obtained graph by $H$. Add
edges to $K_{2, \sigma_2/2}$ such that the part
$\overline{K}_{\sigma_2/2}$  of $K_{2, \sigma_2/2}$ is
$K_{\sigma_2/2}$, we denote the obtained graph by $H'$. Take $m$
copies of $H$, denoted $H_1,\cdots, H_m$ and label the two
non-neighbor vertices of $H_i$ with $x_{i,1}, x_{i,2}$.  Take two
copies of $H'$, denoted $H_0, H_{m+1}$ and similarly label their
vertices. Now, connect $x_{i,2}$ with $x_{i+1,1}$ for $i=0,\cdots,m
$ with an edge. The obtained graph $G$ has $n=(m+2)(\sigma_2/2
+1)+2$ vertices, and $d(x_{i,1})+ d(x_{i,2})=\sigma_2$ for
$i=1,\cdots,m $. It is straightforward to verify that a shortest
path from $x_{0,1}$ to $x_{{m+1},2}$ has length
$3m+5=\frac{6n}{\sigma_2 +2}- \frac{\sigma_2+14}{\sigma_2+2}$

\noindent{\bf Example 2 [3]: } Take $m$ copies of $K_{\delta +1}$,
denoted $X_1,\cdots, X_m$ and label the vertices of $X_i$ with
$x_{i,1},\cdots x_{i,\delta +1}$. Take two copies of $K_{\delta
+2}$, denoted $X_0, X_{m+1}$ and similarly label their vertices.
Now, connect $x_{i,2}$ with $x_{i+1,1}$ for $i=0,\cdots,m $ with an
edge, and delete the edges $(x_{i,1}, x_{i,2})$ for $i=0,\cdots, m+1
$. We can see $d(x_{0,2})+ d(x_{1,2})=\sigma_2=2\delta$. This has
constructed a connected $n$-vertex graph $G$ with $\sigma_2=2\delta
$. The graph has $n=(m+2)(\delta +1)+2$ vertices, and
$diam(G)=3m+5=\frac{6n}{\sigma_2
+2}-(\frac{\sigma_2}{2}+7)/(\frac{\sigma_2}{2}+1)$.

A vertex-colored graph is rainbow vertex-connected if any two
vertices are connected by a path whose internal vertices have
distinct colors. The rainbow vertex-connection of a connected graph
$G$, denoted by $rvc(G)$, is the smallest number of colors that are
needed in order to make $G$ rainbow vertex-connected. The concept of
rainbow vertex-connection was introduced by Krivelevich and Yuster
[7]. It is obvious that if $G$ is a complete graph then $rvc(G) =
0$, if $G$ is a graph of order $n$ then $rvc(G)\leq n-2$.  And
$rvc(G) \geq diam(G)-1$ with equality if the diameter is 1 or 2. In
some case $rvc(G)$ may be much smaller than $rc(G)$. However, in
some other case $rvc(G)$ may be much bigger than $rc(G)$. We may see
some examples given in Krivelevich and Yuster [7] in which they
obtained $rvc(G)\leq \frac{11n}{\delta(G)}$ for every connected
graph with $n$ vertices. Nevertheless, we are able to prove a
theorem analogous to Theorem 1 for the rainbow vertex-connected
case, which is stated as the following Theorem 2. \vskip 0.2 cm
\noindent{\bf Theorem 2.} For a connected graph $G$ of order $n$,
$rvc(G)\leq 8\frac{n-2}{\sigma_2+2}+10$ for $2\leq \sigma_2\leq 6,
\sigma_2\geq 28$, while for $ 7 \leq \sigma_2\leq 8, 16\leq
\sigma_2\leq 27$, $ rvc(G)\leq \frac{10n-16}{\sigma_2+2}+10$, and
for $9 \leq \sigma_2\leq 15,$ $ rvc(G)\leq
\frac{10n-16}{\sigma_2+2}+A(\sigma_2)$ where $ A(\sigma_2)=
63,41,27,20,16,13,11,$ respectively.
 \vskip 0.2 cm

The following notions are needed in the sequel, which could be found
in [5, 7]. Given a graph $G$, a set $D\subseteq V(G)$ is called a
$k$-step dominating set of $G$, if every vertex in $G$ is at a
distance at most $k$ from $D$. Further, if $D$ induces a connected
subgraph of $G$, it is called a connected $k$-step dominating set of
$G$. The $k$-step open neighborhood of a set $D\subseteq V(G)$ is
$N^{k}(D):=\{x\in V(G)| d(x,D)=k \}$, $k=\{0,1,2,\cdots\}$. A
dominating set $D$ in a graph $G$ is called a two-way dominating set
if every pendant vertex of $G $ is included in $D$. In addition, if
$G[D]$ is connected, we call $D$ a connected two-way dominating set.
A connected two-step dominating set $D$ of vertices in a graph $G$
is called a connected two-way two-step dominating set if (i) every
pendant vertex of $G$ is included in $D$ and (ii) every vertex in
$N^{2}(D)$ has at least two neighbors in $N^{1}(D$). We call a
two-step dominating set $k$-strong if every vertex that is not
dominated by it has at least $k$ neighbors that are dominated by it.

\section{ Proof of Theorem 1 }

We start with three lemmas that are needed in order to establish
Theorem 1.

\noindent{\bf Lemma 1.1.} Every connected graph $G$ of order $n$
with at most one pendent vertex has a connected two-step dominating
set $D$ of size at most $6\frac{n-|N^{2}(D)|-2}{\sigma_2+2}+1$ with
equality if $D$ has only one pendent vertex.

\noindent{\bf Proof.} We execute the following stage procedure.

 \noindent{\bf Stage 1.}  $D=\{u\}$, for some $u\in V(G)$ satisfying there exists some vertex  $v\in V(G)$,

  { \quad \quad \quad $uv\not\in E(G), d(u)\geq d(v) $}.

{ \quad \quad \quad While $G[N^{3}(D)]$ is not a complete graph,}

{ \quad \quad \quad \{

{ \quad \quad \quad \quad \quad  pick any $v\in N^{3}(D)$ satisfying
there exists some vertex  $v'\in N^{3}(D)$}

{ \quad \quad \quad  \quad \quad $vv'\not\in E(G),d(v)\geq d(v')$.
Let $(v,v_2,v_1,v_0), v_0\in D$ be a shortest}

{ \quad \quad \quad  \quad \quad   $v-D$ path. $ D=D\cup \{v,v_2,v_1\}$. }

 \quad \quad \quad \} }

Notice that $D$ remains connected after every iteration in Stage 1.
Let $k_1$ be the number of iterations executed in Stage 1. When
Stage 1 starts, $|D\cup N^{1}(D)|\geq \frac{\sigma_2}{2}+1$, since a
new vertex from $N^{3}(D)$ is added to $D$, $|D\cup N^{1}(D)|$
increases by at least $\frac{\sigma_2}{2}+1$ in each iteration, when
Stage 1 ends, $k_1+1\leq \frac{|D\cup
N^{1}(D)|}{\frac{\sigma_2}{2}+1}
=\frac{n-|N^{2}(D)|-|N^{3}(D)|}{\frac{\sigma_2}{2}+1}$.  Since three
more vertices are added in each iteration, $|D|=3k_1+1 \leq
3\frac{n-|N^{2}(D)|-|N^{3}(D)|}{\frac{\sigma_2}{2}+1}-2$.

Initialize $D'=D$, take a vertex $t\in N^{3}(D')$, let
$(t,t_2,t_1,t_0), t_0\in D'$ be a shortest $t-D'$ path,
$D'=D'\cup\{t,t_2,t_1\}$. By this time, $D'$ has been a connected
two-step dominating set. As $|N^{2}(D')|\leq |N^{2}(D)|-1$, if
$|N^{3}(D)|> 1$, then $|D'|<
3\frac{n-|N^{2}(D')|-2}{\frac{\sigma_2}{2}+1}+1$; if $|N^{3}(D)|=
1$, then $|D'|\leq 3\frac{n-|N^{2}(D')|-2}{\frac{\sigma_2}{2}+1}+1$,
and at the same time, we notice that the pendent vertex is in $D'$.
Finally, $D: =D'$, the result follows.  \qed

\noindent{\bf Lemma 1.2.}  Every connected graph $G$ of order $n$
with at most one pendent vertex has a connected two-way two-step
dominating set $D$ of size at most $6\frac{n-2}{\sigma_2+2}+2$.

\noindent{\bf Proof.}  We execute the following stage procedure.

\noindent{\bf Stage 2.}  $D_0=D$ obtained from Stage 1.

 { \quad \quad \quad While $\exists u,v\in N^{2}(D_0),  uv\not\in E(G), d(u, N^{1}(D_0))=d(v, N^{1}(D_0))=1$,}

 { \quad \quad \quad   and $d(u)\geq d(v)$,}

{ \quad \quad \quad \{

{ \quad \quad \quad \quad \quad   $D_0=D_0\cup\{u,u_1\}, (u,u_1,u_0), u_0\in D_0$ be a shortest $u-D_0$ path.}

\quad \quad \quad \} }

Clearly, $D_0$ remains a connected two-step dominating set in Stage
2. Stage 2 ends only when $N^{2}(D_0)$ can be partitioned into two
parts $N_1^{2}(D_0)$ and $N_2^{2}(D_0)$, for any $v\in N_1^{2}(D_0),
d(v, N^{1}(D_0))\geq 2$, and  for any $v\in N_2^{2}(D_0), d(v,
N^{1}(D_0))=1$ and $G[N_2^{2}(D_0)]$ is a complete graph, where $|
N_1^{2}(D_0)|\geq 0, | N_2^{2}(D_0)|\geq 0$.

Let $k_2$ be the number of iterations executed in Stage 2, we add to
$D_0$ a vertex which has at least $\frac{\sigma_2}{2}-1$ neighbors
in $ N^{2}(D_0)$, $ |N^{2}(D_0)|$ reduces by at least
$\frac{\sigma_2}{2}$ in every iteration. Since we start with $|
N^{2}(D)|$ vertices, $k_2\leq \frac{|
N^{2}(D)|}{\frac{\sigma_2}{2}}$. Since we add two vertices to $D_0$
in each iteration, then $|D_0|=|D|+2k_2$, so $|D_0|\leq
6\frac{n-|N^{2}(D)|-2}{\sigma_2+2}+1 +4\frac{| N^{2}(D)|}{\sigma_2}
< 6\frac{n-2}{\sigma_2+2}+1$. We get $|D_0| \leq
6\frac{n-2}{\sigma_2+2}$.

Initialize $D=D_0$, take a vertex $w\in N_2^{2}(D)$, let
$(w,w_1,w_0), w_0\in D$ be a shortest $w-D$ path,
$D=D\cup\{w,w_1\}$, $|D| \leq 6\frac{n-2}{\sigma_2+2}+2$ with the
equality if $D$ has one degree vertex.

If $G$ has no pendent vertex, then $D$ is exactly the two-way
two-step dominating set, so Lemma 1.2 follows. If $G$ has one
pendent vertex, then the pendent vertex is in $ D \cup N^{1}(D)$.
From the above discussion, we know that $D$ is exactly the two-way
two-step dominating set of size at most $6\frac{n-2}{\sigma_2+2}+2$.
\qed

\noindent{\bf Lemma 1.3 [5].} If $D$ is a connected two-way two-step
dominating set in a graph $G$, then $rc(G) \leq rc(G[D]) + 6$.

\noindent{\bf Proof of Theorem 1.} If $G$ has at least two pendent
vertices, then $\sigma_2=2$. As $rc(G)\leq n-1$, however,
$6\frac{n-2}{\sigma_2+2}+8=6\frac{n-2}{2+2}+8 >n-1$, the result is
true. So we may assume that $G$ has at most one pendnet vertex.
Observe that the connected two-way two-step dominating set $D$ can
be rainbow colored using $|D|-1$ colors by ensuring that every edge
of some spanning tree gets distinct colors. So the upper bound
follows immediately from Lemmas 1.2 and 1.3. The tight examples were
given in our introduction. \qed

\section{ Proof of Theorem 2}

\noindent{\bf Lemma 2.1.}  If $G$ is a connected graph of order $n$
with $\sigma_2\geq 12$, then $G$ has a connected
$\frac{\sigma_2}{6}$ -strong two-step dominating set $D$ whose size
is at most $6\frac{n-2}{\sigma_2+2}+2$.

 \noindent{\bf Proof.} We execute the following stage.

 \noindent{\bf Stage 3.}  $D_0=D$ obtained from Stage 1.

 { \quad \quad \quad While $\exists u,v\in N^{2}(D_0),  uv\not\in E(G), d(u,N^{1}(D_0))\leq \frac{\sigma_2}{6}-1$,}

{ \quad \quad \quad   $d(v,N^{1}(D_0))\leq \frac{\sigma_2}{6}-1$ and $d(u)\geq d(v)$,}

{ \quad \quad \quad \{

{ \quad \quad \quad \quad \quad   $ (u,u_1,u_0), u_0\in D_0$ be a shortest $u-D_0$ path, $D_0=D_0\cup\{u,u_1\}$.}

\quad \quad \quad \} }

Notice that $D_0$ remains a connected two-step dominating set in
Stage 3. Stage 3 ends only when $N^{2}(D)$ can be partitioned into
two parts $N_1^{2}(D_0)$ and $ N_2^{2}(D_0)$, for any $v\in
N_1^{2}(D_0)$, $d(v, N^{1}(D_0))\geq \frac{\sigma_2}{6} $, and for
any $v\in N_2^{2}(D_0), d(v, N^{1}(D_0))\leq \frac{\sigma_2}{6}-1$
and $G[N_2^{2}(D_0)]$ is a complete graph, where $|
N_1^{2}(D_0)|\geq 0, | N_2^{2}(D_0)|\geq 0$. Let $k_2$ be the number
of iterations executed in Stage 3, we add to $D_0$ a vertex which
has at least
$\frac{\sigma_2}{2}-\frac{\sigma_2}{6}+1=\frac{\sigma_2}{3}+1$
neighbors in $ N^{2}(D_0)$, $ |N^{2}(D_0)|$ reduces by at least
$\frac{\sigma_2}{3}+2$ in every iteration. We start with $|
N^{2}(D)|$ vertices, so $k_2\leq \frac{|
N^{2}(D)|}{\frac{\sigma_2}{3}+2}$. Since we add two vertices to
$D_0$ in each iteration, then $|D_0|=|D|+2k_2$, so $|D_0|\leq
6\frac{n-|N^{2}(D)|-2}{\sigma_2+2}+1+6\frac{|
N^{2}(D)|}{\sigma_2+6}< 6\frac{n-2}{\sigma_2+2}+1$, hence $|D_0|
\leq 6\frac{n-2}{\sigma_2+2}$.

Initialize $D=D_0$, take a vertex $w\in N_2^{2}(D)$, let
$(w,w_1,w_0), w_0\in D$ be a shortest $w-D$ path,
$D=D\cup\{w,w_1\}$. It is obvious that $D$ also remains connected,
and $|D|\leq 6\frac{n-2}{\sigma_2+2}+2$. \qed

\vskip 0.3 cm

\noindent{\bf Lemma 2.2.} If $G$ is a connected graph of order $n$
with $\sigma_2\geq 9$, then $G$ has a connected $
\frac{\sigma_2}{4}$ -strong two-step dominating set $D$ whose size
is at most $ 8\frac{n-2}{\sigma_2+2}+2$.

\noindent{\bf Proof.} The proof is similar to that of Lemma 2.1, we
execute the following stage 4.

 \noindent{\bf Stage 4.}  $D_0=D$ obtained from Stage 1.

 { \quad \quad \quad While $\exists u,v\in N^{2}(D_0), uv\not\in E(G), d(u,N^{1}(D_0))\leq \frac{\sigma_2}{4}-1$,}

{ \quad \quad \quad   $d(v,N^{1}(D_0))\leq \frac{\sigma_2}{4}-1$ and $d(u)\geq d(v)$,}

{ \quad \quad \quad \{

{ \quad \quad \quad \quad \quad   $ (u,u_1,u_0), u_0\in D_0$ be a shortest $u-D_0$ path, $D_0=D_0\cup\{u,u_1\}$.}

\quad \quad \quad \} }

Notice that $D_0$ remains a connected two-step dominating set in
Stage 4.  Stage 4 ends only when $N^{2}(D)$ can be partitioned into
two parts $N_1^{2}(D_0)$  and $ N_2^{2}(D_0)$, for any $v\in
N_1^{2}(D_0)$, $d(v, N^{1}(D_0))\geq\frac{\sigma_2}{4}$, and for any
$v\in N_2^{2}(D_0), d(v, N^{1}(D_0))\leq \frac{\sigma_2}{4}-1$ and
$G[N_2^{2}(D_0)]$ is a complete graph, where $| N_1^{2}(D_0)|\geq 0,
| N_2^{2}(D_0)|\geq 0$. Let $k_2$ be the number of iterations
executed in Stage 4, we add to $D_0$ a vertex which has at least
$\frac{\sigma_2}{2}-\frac{\sigma_2}{4}+1=\frac{\sigma_2}{4}+1$
neighbors in $ N^{2}(D_0)$, $ |N^{2}(D_0)|$ reduces by at least
$\frac{\sigma_2}{4}+2$ in every iteration. Since we start with $|
N^{2}(D)|$ vertices, $k_2\leq \frac{|
N^{2}(D)|}{\frac{\sigma_2}{4}+2}$. Since we add two vertices to
$D_0$ in each iteration, then $|D_0|=|D|+2k_2$, so $|D_0|\leq
6\frac{n-|N^{2}(D)|-2}{\sigma_2+2}+1+8\frac{|
N^{2}(D)|}{\sigma_2+6}< 6\frac{n-2}{\sigma_2+2}+2\frac{|
N^{2}(D)|}{\sigma_2+2}+1$. As $| N^{2}(D)| \leq n-2$, thus $|D_0|
\leq 8\frac{n-2}{\sigma_2+2}$.

Initialize $D=D_0$, take a vertex $w\in N_2^{2}(D)$, let
$(w,w_1,w_0), w_0\in D$ be a shortest $w-D$ path,
$D=D\cup\{w,w_1\}$. It is obvious that $D$ also remains connected,
and $|D|\leq 8\frac{n-2}{\sigma_2+2}+2$.\qed

\vskip 0.3 cm

\noindent{\bf Lemma 2.3.}  If $G$ is a connected graph of order $n$
with a value of $\sigma_2$, then $G$ has a connected spanning
subgraph with the same value of $\sigma_2$ as $G$ that has less than
$\frac{1}{2}n\sigma_2+\frac{2n}{\sigma_2+4}$ edges.

\noindent{\bf Proof.} If there exist two vertices $u,v\in V(G)$ such
that $uv\not\in E(G), d(u)+d(v)>\sigma_2$, then we delete the edges
incident with the vertices $u,v$ as long as there are any we obtain
a spanning subgraph with $\sigma_2$ and less than
$\frac{1}{2}n\sigma_2 $ edges.  The spanning subgraph has at most
$\frac{n}{\frac{1}{2}\sigma_2+2}$ connected components. Thus  by
adding back at most
$\frac{n}{\frac{1}{2}(\sigma_2+4)}-1=\frac{2n}{\sigma_2+4}-1$ edges,
we can make it connected. \qed

\vskip 0.3 cm

\noindent{\bf Lemma 2.4 (The  Lov\'{a}sz Local Lemma [1]):} Let
$A_1,A_2, \cdots,A_n$ be the events in an arbitrary probability
space. Suppose that each event $A_i$ is mutually independent of a
set of all the other events $A_j$ but at most $d$, and that $P[A_i]
\leq p$ for all $1\leq i\leq n$. If $ep(d + 1) <1$, then
$Pr[\bigwedge^{n}_{i=1}\overline{A_i}]> 0$.

\vskip 0.3 cm \noindent{\bf  Proof of Theorem 2.} Suppose that $G$
is a connected graph with $n$ vertices. By Lemma 2.3 we may assume
that $G$ has less than $\frac{1}{2}n\sigma_2+\frac{2n}{\sigma_2+4}$
edges.

We use Lemma 2.1 to construct a set $D$ which is a  $
\frac{\sigma_2}{6}$ -strong two-step dominating set $D$ whose size
is at most $ 6\frac{n-2}{\sigma_2+2}+2$.

We partition $N^{1}(D)$ into two parts $D_1$ and $D_2$, where $D_1$
are those vertices with at least $\frac{1}{4}(\sigma_2+2)^{2}-1$
neighbors in $N^{2}(D)$. Since $G$ has less than
$\frac{1}{2}n\sigma_2+\frac{2n}{\sigma_2+4}$ edges, we have $|D_1|<
\frac{2n}{\sigma_2+2}$. Denote by $L_1=\{v\in N^{2}(D): v$ has at
least one neighbor in $D_1$\} , and $L_2= N^{2}(D)\setminus L_1$.

We are now ready to describe our coloring. The vertices of $D\cup
D_1$ are each colored with a distinct color. The vertex of $D_2$ are
colored only with 9 fresh colors so that each vertex of $D_2$
chooses its color randomly and independently from all other vertices
of $D_2$. The vertices of $N^{2}(D)$ remain uncolored. Hence, the
total number of colors we used is at most $|D|+|D_1|+9\leq
6\frac{n-2}{\sigma_2+2}+2+\frac{2n}{\sigma_2+2}-1
+9=8\frac{n-2}{\sigma_2+2}+10$.

For each vertex $u$ of $L_2$, let $A_u$ be the event that all the
neighbors of $u $ in $D_2$ are assigned at least two distinct
colors. Now we will prove $Pr[A_u]> 0$ for each $u \in L_2$. Notice
that each vertex $u \in L_2$ has at least $ \frac{\sigma_2}{6}$
neighbors in $D_2$ since $D$ is a connected $\frac{\sigma_2}{6}$
-strong two-step dominating set of $G$. Therefore, we fix a set
$X(u)\subset D_2$ of neighbors of $u$ with
$|X(u)|=\lceil\frac{\sigma_2}{6}\rceil$. Let $B_u$ be the event that
all of the vertices in $X(u)$ receive the same color. Thus,
$Pr[B_u]\leq 9^{-\lceil\frac{\sigma_2}{6}\rceil+1}$. As each vertex
of $D_2$ has less than $\frac{1}{4}(\sigma_2+2)^{2}-1$ neighbors in
$N^{2}(D)$, we have that the event $B_u$ is independent of all other
events $B_v$ for $v\neq u$ but at most $
(\frac{1}{4}(\sigma_2+2)^{2}-2)\lceil\frac{\sigma_2}{6}\rceil$ of
them. Since $e\cdot
9^{-\lceil\frac{\sigma_2}{6}\rceil+1}(((\frac{1}{4}
(\sigma_2+2)^{2}-2)\lceil\frac{\sigma_2}{6}\rceil+1)< 1$ for all
$\sigma_2\geq 28$, by the Lov\'{a}sz Local Lemma, we have $Pr[A_u]>
0$ for each $u \in L_2$. Therefore, for $D_2$, there exists a
coloring with 9 colors such that every vertex of $L_2$ has at least
two neighbors in $D_2$ colored differently.

We now proved the coloring together with the coloring of $D\cup D_1$
with distinct colors, yields a rainbow vertex-connected graph. As
$D\cup D_1$ is connected, and since each vertex of $D_2$ has a
neighbor in $D$, we only need to show that any pair of vertices of
$L$ has a rainbow path connecting them. Notice that each $v\in L$
has at least two neighbors in $N^{1}(D)$ colored differently. Now
let $u,v\in L, x\in N^{1}(D)$ be a neighbor of $u$ and $y\in
N^{1}(D)$ be a neighbor of $v$ whose color is different from the
color of $x$. As there is a rainbow path from $x$ to $y$ whose
internal vertices are only taken from $D$, the result follows.

In the following we still make use of the above $G$, but we use
Lemma 2.2 to construct a set $D$ which is a $\frac{\sigma_2}{4}$
-strong two-step dominating set $D$ whose size is at most $
8\frac{n-2}{\sigma_2+2}+2$.

We still partition $N^{1}(D)$ into two parts $D_1$ and $D_2$, where
$D_1$ are those vertices with at least
$\frac{1}{4}(\sigma_2+2)^{2}-1$ neighbors in  $N^{2}(D)$. We have
$|D_1|< \frac{2n}{\sigma_2+2}$. Denote by $L_1=\{v\in N^{2}(D): v$
has at least one neighbor in $D_1$\}, and $L_2= N^{2}(D)\setminus
L_1$.

Similar to the above coloring, the vertices of $D\cup D_1$ are each
colored with a distinct color. The vertex of $D_2$ are colored only
with 9 fresh colors so that each vertex of $D_2$ chooses its color
randomly and independently from all other vertices of $D_2$. The
vertices of $N^{2}(D)$ remain uncolored. Hence, the total number of
colors we used is at most $|D|+|D_1|+9\leq
8\frac{n-2}{\sigma_2+2}+2+\frac{2n}{\sigma_2+2}-1
+9=\frac{10n-16}{\sigma_2+2}+10$.

For each vertex $u$ of $L_2$, let $A_u$ be the event that all the
neighbors of $u $ in $D_2$ are assigned at least two distinct
colors. Now we will prove $Pr[A_u]> 0$ for each $u \in L_2$. Notice
that each vertex $u \in L_2$ has at least $
\lceil\frac{\sigma_2}{4}\rceil$ neighbors in $D_2$ since $D$ is a
connected $\lceil\frac{\sigma_2}{4}\rceil$ -strong two-step
dominating set of $G$. Therefore, we fix a set $X(u)\subset D_2$ of
neighbors of $ u$ with $|X(u)|=\lceil\frac{\sigma_2}{4}\rceil$. Let
$B_u$ be the event that all of the vertices in $X(u)$ receive the
same color. Thus, $Pr[B_u]\leq
9^{-\lceil\frac{\sigma_2}{4}\rceil+1}$. As each vertex of $D_2$ has
less than $\frac{1}{4}(\sigma_2+2)^{2}-1$ neighbors in $N^{2}(D)$,
we have that the event $B_u$ is independent of all other events
$B_v$ for $v\neq u$ but at most $
(\frac{1}{4}(\sigma_2+2)^{2}-2)\lceil\frac{\sigma_2}{4}\rceil$ of
them. Since
$e\cdot9^{-\lceil\frac{\sigma_2}{4}\rceil+1}(((\frac{1}{4}
(\sigma_2+2)^{2}-2)\lceil\frac{\sigma_2}{4}\rceil+1)< 1$ for all
$\sigma_2\geq 17$, by the Lov\'{a}sz Local Lemma, we have $Pr[A_u]>
0$ for each $u \in L_2$. Therefore, for $D_2$, there exists a
coloring with 9 colors such that every vertex of $L_2$ has at least
two neighbors in $D_2$ colored differently. Similarly, we may show
that $G$ is rainbow vertex-connected. For
$\sigma_2=16,15,14,13,12,11,10,9$ we can use
$10,11,13,16,20,27,41,63 $ colors, respectively, to color $D_2$, and
make $G$ rainbow vertex-connected. The proof of Theorem 2 is now
complete. \qed

\end{document}